
\ifx\shlhetal\undefinedcontrolsequence\let\shlhetal\relax\fi

\input amstex
\expandafter\ifx\csname mathdefs.tex\endcsname\relax
  \expandafter\gdef\csname mathdefs.tex\endcsname{}
\else \message{Hey!  Apparently you were trying to
  \string\input{mathdefs.tex} twice.   This does not make sense.} 
\errmessage{Please edit your file (probably \jobname.tex) and remove
any duplicate ``\string\input'' lines}\endinput\fi




\catcode`\X=12\catcode`\@=11

\def\n@wcount{\alloc@0\count\countdef\insc@unt}
\def\n@wwrite{\alloc@7\write\chardef\sixt@@n}
\def\n@wread{\alloc@6\read\chardef\sixt@@n}
\def\r@s@t{\relax}\def\v@idline{\par}\def\@mputate#1/{#1}
\def\l@c@l#1X{\firstpart.#1}\def\gl@b@l#1X{#1}\def\t@d@l#1X{{}}

\def\crossrefs#1{\ifx\all#1\let\tr@ce=\all\else\def\tr@ce{#1,}\fi
   \n@wwrite\cit@tionsout\openout\cit@tionsout=\jobname.cit 
   \write\cit@tionsout{\tr@ce}\expandafter\setfl@gs\tr@ce,}
\def\setfl@gs#1,{\def\@{#1}\ifx\@\empty\let\next=\relax
   \else\let\next=\setfl@gs\expandafter\xdef
   \csname#1tr@cetrue\endcsname{}\fi\next}
\def\m@ketag#1#2{\expandafter\n@wcount\csname#2tagno\endcsname
     \csname#2tagno\endcsname=0\let\tail=\all\xdef\all{\tail#2,}
   \ifx#1\l@c@l\let\tail=\r@s@t\xdef\r@s@t{\csname#2tagno\endcsname=0\tail}\fi
   \expandafter\gdef\csname#2cite\endcsname##1{\expandafter
     \ifx\csname#2tag##1\endcsname\relax?\else\csname#2tag##1\endcsname\fi
     \expandafter\ifx\csname#2tr@cetrue\endcsname\relax\else
     \write\cit@tionsout{#2tag ##1 cited on page \folio.}\fi}
   \expandafter\gdef\csname#2page\endcsname##1{\expandafter
     \ifx\csname#2page##1\endcsname\relax?\else\csname#2page##1\endcsname\fi
     \expandafter\ifx\csname#2tr@cetrue\endcsname\relax\else
     \write\cit@tionsout{#2tag ##1 cited on page \folio.}\fi}
   \expandafter\gdef\csname#2tag\endcsname##1{\expandafter
      \ifx\csname#2check##1\endcsname\relax
      \expandafter\xdef\csname#2check##1\endcsname{}%
      \else\immediate\write16{Warning: #2tag ##1 used more than once.}\fi
      \multit@g{#1}{#2}##1/X%
      \write\t@gsout{#2tag ##1 assigned number \csname#2tag##1\endcsname\space
      on page \number\count0.}%
   \csname#2tag##1\endcsname}}

\def\multit@g#1#2#3/#4X{\def\t@mp{#4}\ifx\t@mp\empty%
      \global\advance\csname#2tagno\endcsname by 1 
      \expandafter\xdef\csname#2tag#3\endcsname
      {#1\number\csname#2tagno\endcsnameX}%
   \else\expandafter\ifx\csname#2last#3\endcsname\relax
      \expandafter\n@wcount\csname#2last#3\endcsname
      \global\advance\csname#2tagno\endcsname by 1 
      \expandafter\xdef\csname#2tag#3\endcsname
      {#1\number\csname#2tagno\endcsnameX}
      \write\t@gsout{#2tag #3 assigned number \csname#2tag#3\endcsname\space
      on page \number\count0.}\fi
   \global\advance\csname#2last#3\endcsname by 1
   \def\t@mp{\expandafter\xdef\csname#2tag#3/}%
   \expandafter\t@mp\@mputate#4\endcsname
   {\csname#2tag#3\endcsname\lastpart{\csname#2last#3\endcsname}}\fi}
\def\t@gs#1{\def\all{}\m@ketag#1e\m@ketag#1s\m@ketag\t@d@l p
\let\realscite\scite
\let\realstag\stag
   \m@ketag\gl@b@l r \n@wread\t@gsin
   \openin\t@gsin=\jobname.tgs \re@der \closein\t@gsin
   \n@wwrite\t@gsout\openout\t@gsout=\jobname.tgs }
\outer\def\localtags{\t@gs\l@c@l}
\outer\def\globaltags{\t@gs\gl@b@l}
\outer\def\newlocaltag#1{\m@ketag\l@c@l{#1}}
\outer\def\newglobaltag#1{\m@ketag\gl@b@l{#1}}

\newif\ifpr@ 
\def\m@kecs #1tag #2 assigned number #3 on page #4.%
   {\expandafter\gdef\csname#1tag#2\endcsname{#3}
   \expandafter\gdef\csname#1page#2\endcsname{#4}
   \ifpr@\expandafter\xdef\csname#1check#2\endcsname{}\fi}
\def\re@der{\ifeof\t@gsin\let\next=\relax\else
   \read\t@gsin to\t@gline\ifx\t@gline\v@idline\else
   \expandafter\m@kecs \t@gline\fi\let \next=\re@der\fi\next}
\def\pretags#1{\pr@true\pret@gs#1,,}
\def\pret@gs#1,{\def\@{#1}\ifx\@\empty\let\n@xtfile=\relax
   \else\let\n@xtfile=\pret@gs \openin\t@gsin=#1.tgs \message{#1} \re@der 
   \closein\t@gsin\fi \n@xtfile}

\newcount\sectno\sectno=0\newcount\subsectno\subsectno=0
\newif\ifultr@local \def\ultralocal{\ultr@localtrue}
\def\firstpart{\number\sectno}
\def\lastpart#1{\ifcase#1 \or a\or b\or c\or d\or e\or f\or g\or h\or 
   i\or k\or l\or m\or n\or o\or p\or q\or r\or s\or t\or u\or v\or w\or 
   x\or y\or z \fi}

\def\resetall{\global\advance\sectno by 1\subsectno=0
   \gdef\firstpart{\number\sectno}\r@s@t}
\def\resetsub{\global\advance\subsectno by 1
   \gdef\firstpart{\number\sectno.\number\subsectno}\r@s@t}
\def\newsection#1\par{\resetall\vskip0pt plus.3\vsize\penalty-250
   \vskip0pt plus-.3\vsize\bigskip\bigskip
   \message{#1}\leftline{\bf#1}\nobreak\bigskip}
\def\subsection#1\par{\ifultr@local\resetsub\fi
   \vskip0pt plus.2\vsize\penalty-250\vskip0pt plus-.2\vsize
   \bigskip\smallskip\message{#1}\leftline{\bf#1}\nobreak\medskip}


\newdimen\marginshift

\newdimen\margindelta
\newdimen\marginmax
\newdimen\marginmin

\def\margininit{       
\marginmax=3 true cm                  
				      
\margindelta=0.1 true cm              
\marginmin=0.1true cm                 
\marginshift=\marginmin
}    

\def\t@gsjj#1,{\def\@{#1}\ifx\@\empty\let\next=\relax\else\let\next=\t@gsjj
   \def\@@{p}\ifx\@\@@\else
   \expandafter\gdef\csname#1cite\endcsname##1{\citejj{##1}}
   \expandafter\gdef\csname#1page\endcsname##1{?}
   \expandafter\gdef\csname#1tag\endcsname##1{\tagjj{##1}}\fi\fi\next}
\newif\ifshowstuffinmargin
\showstuffinmarginfalse
\def\jjtags{\ifx\shlhetal\relax 
  \else
\ifx\shlhetal\undefinedcontrolseq
\else
\showstuffinmargintrue
\ifx\all\relax\else\expandafter\t@gsjj\all,\fi\fi \fi
}

\def\tagjj#1{\realstag{#1}\oldmginpar{\zeigen{#1}}}
\def\citejj#1{\rechnen{#1}\mginpar{\zeigen{#1}}}     

\def\rechnen#1{\expandafter\ifx\csname stag#1\endcsname\relax ??\else
                           \csname stag#1\endcsname\fi}

\newdimen\theight

\def\marginfont{\sevenrm}

\def\trymarginbox#1{\setbox0=\hbox{\marginfont\hskip\marginshift #1}%
		\global\marginshift\wd0 
		\global\advance\marginshift\margindelta}

\def \oldmginpar#1{%
\ifvmode\setbox0\hbox to \hsize{\hfill\rlap{\marginfont\quad#1}}%
\ht0 0cm
\dp0 0cm
\box0\vskip-\baselineskip
\else 
             \vadjust{\trymarginbox{#1}%
		\ifdim\marginshift>\marginmax \global\marginshift\marginmin
			\trymarginbox{#1}%
                \fi
             \theight=\ht0
             \advance\theight by \dp0    \advance\theight by \lineskip
             \kern -\theight \vbox to \theight{\rightline{\rlap{\box0}}%
\vss}}\fi}

\newdimen\upordown
\global\upordown=8pt
\font\tinyfont=cmtt8 
\def\mginpar#1{\smash{\hbox to 0cm{\kern-10pt\raise7pt\hbox{\tinyfont #1}\hss}}}
\def\mginpar#1{{\hbox to 0cm{\kern-10pt\raise\upordown\hbox{\tinyfont #1}\hss}}\global\upordown-\upordown}


\def\t@gsoff#1,{\def\@{#1}\ifx\@\empty\let\next=\relax\else\let\next=\t@gsoff
   \def\@@{p}\ifx\@\@@\else
   \expandafter\gdef\csname#1cite\endcsname##1{\zeigen{##1}}
   \expandafter\gdef\csname#1page\endcsname##1{?}
   \expandafter\gdef\csname#1tag\endcsname##1{\zeigen{##1}}\fi\fi\next}
\def\verbatimtags{\showstuffinmarginfalse
\ifx\all\relax\else\expandafter\t@gsoff\all,\fi}
\def\zeigen#1{\hbox{$\scriptstyle\langle$}#1\hbox{$\scriptstyle\rangle$}}


\def\margintag#1{\ifshowstuffinmargin\oldmginpar{\zeigen{#1}}\fi}

\def\(#1){\edef\dot@g{\ifmmode\ifinner(\hbox{\noexpand\etag{#1}})
   \else\noexpand\eqno(\hbox{\noexpand\etag{#1}})\fi
   \else(\noexpand\ecite{#1})\fi}\dot@g}

\newif\ifbr@ck
\def\eat#1{}
\def\[#1]{\br@cktrue[\br@cket#1'X]}
\def\br@cket#1'#2X{\def\temp{#2}\ifx\temp\empty\let\next\eat
   \else\let\next\br@cket\fi
   \ifbr@ck\br@ckfalse\br@ck@t#1,X\else\br@cktrue#1\fi\next#2X}
\def\br@ck@t#1,#2X{\def\temp{#2}\ifx\temp\empty\let\neext\eat
   \else\let\neext\br@ck@t\def\temp{,}\fi
   \def\teemp{#1}\ifx\teemp\empty\else\rcite{#1}\fi\temp\neext#2X}
\def\resetbr@cket{\gdef\[##1]{[\rtag{##1}]}}
\def\references{\resetbr@cket\newsection References\par}

\newtoks\symb@ls\newtoks\s@mb@ls\newtoks\p@gelist\n@wcount\ftn@mber
    \ftn@mber=1\newif\ifftn@mbers\ftn@mbersfalse\newif\ifbyp@ge\byp@gefalse
\def\defm@rk{\ifftn@mbers\n@mberm@rk\else\symb@lm@rk\fi}
\def\n@mberm@rk{\xdef\m@rk{{\the\ftn@mber}}%
    \global\advance\ftn@mber by 1 }
\def\rot@te#1{\let\temp=#1\global#1=\expandafter\r@t@te\the\temp,X}
\def\r@t@te#1,#2X{{#2#1}\xdef\m@rk{{#1}}}
\def\b@@st#1{{$^{#1}$}}\def\str@p#1{#1}
\def\symb@lm@rk{\ifbyp@ge\rot@te\p@gelist\ifnum\expandafter\str@p\m@rk=1 
    \s@mb@ls=\symb@ls\fi\write\f@nsout{\number\count0}\fi \rot@te\s@mb@ls}
\def\byp@ge{\byp@getrue\n@wwrite\f@nsin\openin\f@nsin=\jobname.fns 
    \n@wcount\currentp@ge\currentp@ge=0\p@gelist={0}
    \re@dfns\closein\f@nsin\rot@te\p@gelist
    \n@wread\f@nsout\openout\f@nsout=\jobname.fns }
\def\m@kelist#1X#2{{#1,#2}}
\def\re@dfns{\ifeof\f@nsin\let\next=\relax\else\read\f@nsin to \f@nline
    \ifx\f@nline\v@idline\else\let\t@mplist=\p@gelist
    \ifnum\currentp@ge=\f@nline
    \global\p@gelist=\expandafter\m@kelist\the\t@mplistX0
    \else\currentp@ge=\f@nline
    \global\p@gelist=\expandafter\m@kelist\the\t@mplistX1\fi\fi
    \let\next=\re@dfns\fi\next}
\def\symbols#1{\symb@ls={#1}\s@mb@ls=\symb@ls} 
\def\bigsymbol{\textstyle}
\symbols{\bigsymbol\ast,\dagger,\ddagger,\sharp,\flat,\natural,\star}
\def\ftnumbers{\ftn@mberstrue} \def\ftsymbols{\ftn@mbersfalse}
\def\paginal{\byp@ge} \def\resetftnumbers{\ftn@mber=1}
\def\ftnote#1{\defm@rk\expandafter\expandafter\expandafter\footnote
    \expandafter\b@@st\m@rk{#1}}

\long\def\jump#1\endjump{}
\def\ssum{\mathop{\lower .1em\hbox{$\textstyle\Sigma$}}\nolimits}

\def\qed{\nobreak\kern 1em \vrule height .5em width .5em depth 0em}
\def\newneq{\hbox{\rlap{\hbox to 1\wd9{\hss$=$\hss}}\raise .1em 
   \hbox to 1\wd9{\hss$\scriptscriptstyle/$\hss}}}
\def\subsetne{\setbox9 = \hbox{$\subset$}\mathrel{\hbox{\rlap
   {\lower .4em \newneq}\raise .13em \hbox{$\subset$}}}}
\def\supsetne{\setbox9 = \hbox{$\subset$}\mathrel{\hbox{\rlap
   {\lower .4em \newneq}\raise .13em \hbox{$\supset$}}}}

\def\vbar{\mathchoice{\vrule height6.3ptdepth-.5ptwidth.8pt\kern-.8pt}
   {\vrule height6.3ptdepth-.5ptwidth.8pt\kern-.8pt}
   {\vrule height4.1ptdepth-.35ptwidth.6pt\kern-.6pt}
   {\vrule height3.1ptdepth-.25ptwidth.5pt\kern-.5pt}}
\def\f@dge{\mathchoice{}{}{\mkern.5mu}{\mkern.8mu}}
\def\b@c#1#2{{\rm \mkern#2mu\vbar\mkern-#2mu#1}}
\def\b@b#1{{\rm I\mkern-3.5mu #1}}
\def\b@a#1#2{{\rm #1\mkern-#2mu\f@dge #1}}
\def\bb#1{{\count4=`#1 \advance\count4by-64 \ifcase\count4\or\b@a A{11.5}\or
   \b@b B\or\b@c C{5}\or\b@b D\or\b@b E\or\b@b F \or\b@c G{5}\or\b@b H\or
   \b@b I\or\b@c J{3}\or\b@b K\or\b@b L \or\b@b M\or\b@b N\or\b@c O{5} \or
   \b@b P\or\b@c Q{5}\or\b@b R\or\b@a S{8}\or\b@a T{10.5}\or\b@c U{5}\or
   \b@a V{12}\or\b@a W{16.5}\or\b@a X{11}\or\b@a Y{11.7}\or\b@a Z{7.5}\fi}}

\catcode`\X=11 \catcode`\@=12




\let\thischap\jobname

\def\partof#1{\csname returnthe#1part\endcsname}
\def\CHAPOF#1{\csname returnthe#1chap\endcsname}

\def\chapof#1{\CHAPOF{#1}}

\def\setchapter#1,#2,#3;{%
  \expandafter\def\csname returnthe#1part\endcsname{#2}%
  \expandafter\def\csname returnthe#1chap\endcsname{#3}%
}

\def\setprevious#1 #2 {%
  \expandafter\def\csname set#1page\endcsname{\input page-#2}
}


 \setchapter  E53,B,N;       \setprevious E53 null
 \setchapter  300z,B,A;       \setprevious 300z E53
 \setchapter  88r,B,I;       \setprevious 88r 300z
 \setchapter  600,B,II;       \setprevious  600 88r
 \setchapter  705,B,III;       \setprevious   705 600
 \setchapter  734,B,IV;        \setprevious   734 705
 \setchapter  300x,B,;      \setprevious   300x 734

 \setchapter 300a,A,V.A;      \setprevious 300a 88r
 \setchapter 300b,A,V.B;       \setprevious 300b 300a
 \setchapter 300c,A,V.C;       \setprevious 300c 300b
 \setchapter 300d,A,V.D;       \setprevious 300d 300c
 \setchapter 300e,A,V.E;       \setprevious 300e 300d
 \setchapter 300f,A,V.F;       \setprevious 300f 300e
 \setchapter 300g,A,V.G;       \setprevious 300g 300f

  \setchapter  E46,B,VI;      \setprevious    E46 734
  \setchapter  838,B,VII;      \setprevious   838 E46

\newwrite\pageout
\def\rememberpagenumber{\let\setpage\relax
\openout\pageout page-\jobname  \relax \write\pageout{\setpage\the\pageno.}}

\def\recallpagenumber{\csname set\jobname page\endcsname
\def\headmark##1{\rightheadtext{\chapof{\jobname}.##1}}\WRITETOC}
\def\setupchapter#1{\leftheadtext{\chapof{\jobname}. #1}}

\def\setpage#1.{\pageno=#1\relax\advance\pageno1\relax}

\def\cprefix#1{
\edef\theotherpart{\partof{#1}}\edef\theotherchap{\chapof{#1}}%
\ifx\theotherpart\thispart
   \ifx\theotherchap\thischap 
    \else 
     \theotherchap%
    \fi
   \else 
     \theotherchap\fi}

\def\sectioncite[#1]#2{%
     \cprefix{#2}#1}

\edef\thispart{\partof{\thischap}}
\edef\thischap{\chapof{\thischap}}

\def\lastpage of '#1' is #2.{\expandafter\def\csname lastpage#1\endcsname{#2}}

\def\yCITE[#1]#2{\cprefix{#2}.\scite{#2-#1}}

\newwrite\writetoc
\immediate\openout\writetoc \jobname.toc
\def\addcontents#1{\def\WRITETOC{\immediate\write\writetoc{\noexpand\tocentry{\chapof{\jobname}}{#1}{\number\pageno}}}}



\def\spuriousreset{}


\expandafter\ifx\csname citeadd.tex\endcsname\relax
\expandafter\gdef\csname citeadd.tex\endcsname{}
\else \message{Hey!  Apparently you were trying to
\string\input{citeadd.tex} twice.   This does not make sense.} 
\errmessage{Please edit your file (probably \jobname.tex) and remove
any duplicate ``\string\input'' lines}\endinput\fi

\sectno=-1   
\localtags
\jjtags
\NoBlackBoxes
\define\mr{\medskip\roster}
\define\sn{\smallskip\noindent}
\define\mn{\medskip\noindent}
\define\bn{\bigskip\noindent}
\define\ub{\underbar}
\define\wilog{\text{without loss of generality}}
\define\ermn{\endroster\medskip\noindent}

\define\dbcu{\dsize\bigcup}
\define\rest{\restriction}

\define \nl{\newline}
\magnification=\magstep 1
\documentstyle{amsppt}

{    
\catcode`@11

\ifx\alicetwothousandloaded@\relax
  \endinput\else\global\let\alicetwothousandloaded@\relax\fi

\gdef\subjclass{\let\savedef@\subjclass
 \def\subjclass##1\endsubjclass{\let\subjclass\savedef@
   \toks@{\def\usualspace{{\rm\enspace}}\eightpoint}%
   \toks@@{##1\unskip.}%
   \edef\thesubjclass@{\the\toks@
     \frills@{{\noexpand\rm2000 {\noexpand\it Mathematics Subject
       Classification}.\noexpand\enspace}}%
     \the\toks@@}}%
  \nofrillscheck\subjclass}
} 


\expandafter\ifx\csname alice2jlem.tex\endcsname\relax
  \expandafter\xdef\csname alice2jlem.tex\endcsname{\the\catcode`@}
\else \message{Hey!  Apparently you were trying to
\string\input{alice2jlem.tex}  twice.   This does not make sense.}
\errmessage{Please edit your file (probably \jobname.tex) and remove
any duplicate ``\string\input'' lines}\endinput\fi

\expandafter\ifx\csname bib4plain.tex\endcsname\relax
  \expandafter\gdef\csname bib4plain.tex\endcsname{}
\else \message{Hey!  Apparently you were trying to \string\input
  bib4plain.tex twice.   This does not make sense.}
\errmessage{Please edit your file (probably \jobname.tex) and remove
any duplicate ``\string\input'' lines}\endinput\fi

\def\renewcommand{\newcommand}	       
\edef\cite{\the\catcode`@}%
\catcode`@ = 11
\let\@oldatcatcode = \cite
\chardef\@letter = 11
\chardef\@other = 12
%
%
%
%
\def\@innerdef#1#2{\edef#1{\expandafter\noexpand\csname #2\endcsname}}%
%
%
\@innerdef\@innernewcount{newcount}%
\@innerdef\@innernewdimen{newdimen}%
\@innerdef\@innernewif{newif}%
\@innerdef\@innernewwrite{newwrite}%
%
%
%
\def\@gobble#1{}%
%
%
%
\ifx\inputlineno\@undefined
   \let\@linenumber = \empty 
\else
   \def\@linenumber{\the\inputlineno:\space}%
\fi
%
%
%
\def\@futurenonspacelet#1{\def\cs{#1}%
   \afterassignment\@stepone\let\@nexttoken=
}%
\begingroup 
\def\\{\global\let\@stoken= }%
\\ 
\endgroup
\def\@stepone{\expandafter\futurelet\cs\@steptwo}%
\def\@steptwo{\expandafter\ifx\cs\@stoken\let\@@next=\@stepthree
   \else\let\@@next=\@nexttoken\fi \@@next}%
\def\@stepthree{\afterassignment\@stepone\let\@@next= }%
%
%
%
\def\@getoptionalarg#1{%
   \let\@optionaltemp = #1%
   \let\@optionalnext = \relax
   \@futurenonspacelet\@optionalnext\@bracketcheck
}%
%
%
\def\@bracketcheck{%
   \ifx [\@optionalnext
      \expandafter\@@getoptionalarg
   \else
      \let\@optionalarg = \empty
      \expandafter\@optionaltemp
   \fi
}%
\def\@@getoptionalarg[#1]{%
   \def\@optionalarg{#1}%
   \@optionaltemp
}%
%
%
%
\def\@nnil{\@nil}%
\def\@fornoop#1\@@#2#3{}%
\def\@for#1:=#2\do#3{%
   \edef\@fortmp{#2}%
   \ifx\@fortmp\empty \else
      \expandafter\@forloop#2,\@nil,\@nil\@@#1{#3}%
   \fi
}%
\def\@forloop#1,#2,#3\@@#4#5{\def#4{#1}\ifx #4\@nnil \else
       #5\def#4{#2}\ifx #4\@nnil \else#5\@iforloop #3\@@#4{#5}\fi\fi
}%
\def\@iforloop#1,#2\@@#3#4{\def#3{#1}\ifx #3\@nnil
       \let\@nextwhile=\@fornoop \else
      #4\relax\let\@nextwhile=\@iforloop\fi\@nextwhile#2\@@#3{#4}%
}%
%
%
%
\@innernewif\if@fileexists
\def\@testfileexistence{\@getoptionalarg\@finishtestfileexistence}%
\def\@finishtestfileexistence#1{%
   \begingroup
      \def\extension{#1}%
      \immediate\openin0 =
         \ifx\@optionalarg\empty\jobname\else\@optionalarg\fi
         \ifx\extension\empty \else .#1\fi
         \space
      \ifeof 0
         \global\@fileexistsfalse
      \else
         \global\@fileexiststrue
      \fi
      \immediate\closein0
   \endgroup
}%
%
%
%
%
\def\bibliographystyle#1{%
   \@readauxfile
   \@writeaux{\string\bibstyle{#1}}%
}%
\let\bibstyle = \@gobble
%
%
\let\bblfilebasename = \jobname
\def\bibliography#1{%
   \@readauxfile
   \@writeaux{\string\bibdata{#1}}%
   \@testfileexistence[\bblfilebasename]{bbl}%
   \if@fileexists
      \nobreak
      \@readbblfile
   \fi
}%
\let\bibdata = \@gobble
%
%
\def\nocite#1{%
   \@readauxfile
   \@writeaux{\string\citation{#1}}%
}%
\@innernewif\if@notfirstcitation
%
%
\def\cite{\@getoptionalarg\@cite}%
%
%
\def\@cite#1{%
   \let\@citenotetext = \@optionalarg
   \printcitestart
   \nocite{#1}%
   \@notfirstcitationfalse
   \@for \@citation :=#1\do
   {%
      \expandafter\@onecitation\@citation\@@
   }%
   \ifx\empty\@citenotetext\else
      \printcitenote{\@citenotetext}%
   \fi
   \printcitefinish
}%
\newif\ifweareinprivate
\weareinprivatetrue
\ifx\shlhetal\undefinedcontrolseq\weareinprivatefalse\fi
\ifx\shlhetal\relax\weareinprivatefalse\fi
\def\@onecitation#1\@@{%
   \if@notfirstcitation
      \printbetweencitations
   \fi
   \expandafter \ifx \csname\@citelabel{#1}\endcsname \relax
      \if@citewarning
         \message{\@linenumber Undefined citation `#1'.}%
      \fi
     \ifweareinprivate
      \expandafter\gdef\csname\@citelabel{#1}\endcsname{%
\strut 
\vadjust{\vskip-\dp\strutbox
\vbox to 0pt{\vss\parindent0cm \leftskip=\hsize 
\advance\leftskip3mm
\advance\hsize 4cm\strut\openup-4pt 
\rightskip 0cm plus 1cm minus 0.5cm ?  #1 ?\strut}}
         {\tt
            \escapechar = -1
            \nobreak\hskip0pt\pfeilsw
            \expandafter\string\csname#1\endcsname
             \pfeilso
            \nobreak\hskip0pt
         }%
      }%
     \else  
      \expandafter\gdef\csname\@citelabel{#1}\endcsname{%
            {\tt\expandafter\string\csname#1\endcsname}
      }%
     \fi  
   \fi
   \csname\@citelabel{#1}\endcsname
   \@notfirstcitationtrue
}%
%
%
\def\@citelabel#1{b@#1}%
%
%
\def\@citedef#1#2{\expandafter\gdef\csname\@citelabel{#1}\endcsname{#2}}%
%
%
%
\def\@readbblfile{%
   \ifx\@itemnum\@undefined
      \@innernewcount\@itemnum
   \fi
   \begingroup
      \def\begin##1##2{%
         \setbox0 = \hbox{\biblabelcontents{##2}}%
         \biblabelwidth = \wd0
      }%
      \def\end##1{}
      %
      %
      \@itemnum = 0
      \def\bibitem{\@getoptionalarg\@bibitem}%
      \def\@bibitem{%
         \ifx\@optionalarg\empty
            \expandafter\@numberedbibitem
         \else
            \expandafter\@alphabibitem
         \fi
      }%
      \def\@alphabibitem##1{%
         \expandafter \xdef\csname\@citelabel{##1}\endcsname {\@optionalarg}%
         \ifx\biblabelprecontents\@undefined
            \let\biblabelprecontents = \relax
         \fi
         \ifx\biblabelpostcontents\@undefined
            \let\biblabelpostcontents = \hss
         \fi
         \@finishbibitem{##1}%
      }%
      \def\@numberedbibitem##1{%
         \advance\@itemnum by 1
         \expandafter \xdef\csname\@citelabel{##1}\endcsname{\number\@itemnum}%
         \ifx\biblabelprecontents\@undefined
            \let\biblabelprecontents = \hss
         \fi
         \ifx\biblabelpostcontents\@undefined
            \let\biblabelpostcontents = \relax
         \fi
         \@finishbibitem{##1}%
      }%
      \def\@finishbibitem##1{%
         \biblabelprint{\csname\@citelabel{##1}\endcsname}%
         \@writeaux{\string\@citedef{##1}{\csname\@citelabel{##1}\endcsname}}%
         \ignorespaces
      }%
      %
      %
      \let\em = \bblem
      \let\newblock = \bblnewblock
      \let\sc = \bblsc
      \frenchspacing
      \clubpenalty = 4000 \widowpenalty = 4000
      \tolerance = 10000 \hfuzz = .5pt
      \everypar = {\hangindent = \biblabelwidth
                      \advance\hangindent by \biblabelextraspace}%
      \bblrm
      \parskip = 1.5ex plus .5ex minus .5ex
      \biblabelextraspace = .5em
      \bblhook
      \input \bblfilebasename.bbl
   \endgroup
}%
%
%
\@innernewdimen\biblabelwidth
\@innernewdimen\biblabelextraspace
%
%
%
\def\biblabelprint#1{%
   \noindent
   \hbox to \biblabelwidth{%
      \biblabelprecontents
      \biblabelcontents{#1}%
      \biblabelpostcontents
   }%
   \kern\biblabelextraspace
}%
%
%
%
\def\biblabelcontents#1{{\bblrm [#1]}}%
%
%
\def\bblrm{\rm}%
%
%
\def\bblem{\it}%
%
%
\def\bblsc{\ifx\@scfont\@undefined
              \font\@scfont = cmcsc10
           \fi
           \@scfont
}%
%
%
\def\bblnewblock{\hskip .11em plus .33em minus .07em }%
%
%
\let\bblhook = \empty
%
%
%
\def\printcitestart{[}
\def\printcitefinish{]}
\def\printbetweencitations{, }
\def\printcitenote#1{, #1}
%
%
%
\let\citation = \@gobble
%
%
%
\@innernewcount\@numparams
%
%
\def\newcommand#1{%
   \def\@commandname{#1}%
   \@getoptionalarg\@continuenewcommand
}%
%
%
\def\@continuenewcommand{%
   \@numparams = \ifx\@optionalarg\empty 0\else\@optionalarg \fi \relax
   \@newcommand
}%
%
%
\def\@newcommand#1{%
   \def\@startdef{\expandafter\edef\@commandname}%
   \ifnum\@numparams=0
      \let\@paramdef = \empty
   \else
      \ifnum\@numparams>9
         \errmessage{\the\@numparams\space is too many parameters}%
      \else
         \ifnum\@numparams<0
            \errmessage{\the\@numparams\space is too few parameters}%
         \else
            \edef\@paramdef{%
               \ifcase\@numparams
                  \empty  No arguments.
               \or ####1%
               \or ####1####2%
               \or ####1####2####3%
               \or ####1####2####3####4%
               \or ####1####2####3####4####5%
               \or ####1####2####3####4####5####6%
               \or ####1####2####3####4####5####6####7%
               \or ####1####2####3####4####5####6####7####8%
               \or ####1####2####3####4####5####6####7####8####9%
               \fi
            }%
         \fi
      \fi
   \fi
   \expandafter\@startdef\@paramdef{#1}%
}%
%
%
%
%
\def\@readauxfile{%
   \if@auxfiledone \else 
      \global\@auxfiledonetrue
      \@testfileexistence{aux}%
      \if@fileexists
         \begingroup
            \endlinechar = -1
            \catcode`@ = 11
            \input \jobname.aux
         \endgroup
      \else
         \message{\@undefinedmessage}%
         \global\@citewarningfalse
      \fi
      \immediate\openout\@auxfile = \jobname.aux
   \fi
}%
%
%
\newif\if@auxfiledone
\ifx\noauxfile\@undefined \else \@auxfiledonetrue\fi
%
%
%
%
\@innernewwrite\@auxfile
\def\@writeaux#1{\ifx\noauxfile\@undefined \write\@auxfile{#1}\fi}%
%
%
%
\ifx\@undefinedmessage\@undefined
   \def\@undefinedmessage{No .aux file; I won't give you warnings about
                          undefined citations.}%
\fi
%
%
\@innernewif\if@citewarning
\ifx\noauxfile\@undefined \@citewarningtrue\fi
%
%
%
\catcode`@ = \@oldatcatcode

\def\pfeilso{\leavevmode
            \vrule width 1pt height9pt depth 0pt\relax
           \vrule width 1pt height8.7pt depth 0pt\relax
           \vrule width 1pt height8.3pt depth 0pt\relax
           \vrule width 1pt height8.0pt depth 0pt\relax
           \vrule width 1pt height7.7pt depth 0pt\relax
            \vrule width 1pt height7.3pt depth 0pt\relax
            \vrule width 1pt height7.0pt depth 0pt\relax
            \vrule width 1pt height6.7pt depth 0pt\relax
            \vrule width 1pt height6.3pt depth 0pt\relax
            \vrule width 1pt height6.0pt depth 0pt\relax
            \vrule width 1pt height5.7pt depth 0pt\relax
            \vrule width 1pt height5.3pt depth 0pt\relax
            \vrule width 1pt height5.0pt depth 0pt\relax
            \vrule width 1pt height4.7pt depth 0pt\relax
            \vrule width 1pt height4.3pt depth 0pt\relax
            \vrule width 1pt height4.0pt depth 0pt\relax
            \vrule width 1pt height3.7pt depth 0pt\relax
            \vrule width 1pt height3.3pt depth 0pt\relax
            \vrule width 1pt height3.0pt depth 0pt\relax
            \vrule width 1pt height2.7pt depth 0pt\relax
            \vrule width 1pt height2.3pt depth 0pt\relax
            \vrule width 1pt height2.0pt depth 0pt\relax
            \vrule width 1pt height1.7pt depth 0pt\relax
            \vrule width 1pt height1.3pt depth 0pt\relax
            \vrule width 1pt height1.0pt depth 0pt\relax
            \vrule width 1pt height0.7pt depth 0pt\relax
            \vrule width 1pt height0.3pt depth 0pt\relax}

\def\pfeilsw{ \leavevmode 
            \vrule width 1pt height0.3pt depth 0pt\relax
            \vrule width 1pt height0.7pt depth 0pt\relax
            \vrule width 1pt height1.0pt depth 0pt\relax
            \vrule width 1pt height1.3pt depth 0pt\relax
            \vrule width 1pt height1.7pt depth 0pt\relax
            \vrule width 1pt height2.0pt depth 0pt\relax
            \vrule width 1pt height2.3pt depth 0pt\relax
            \vrule width 1pt height2.7pt depth 0pt\relax
            \vrule width 1pt height3.0pt depth 0pt\relax
            \vrule width 1pt height3.3pt depth 0pt\relax
            \vrule width 1pt height3.7pt depth 0pt\relax
            \vrule width 1pt height4.0pt depth 0pt\relax
            \vrule width 1pt height4.3pt depth 0pt\relax
            \vrule width 1pt height4.7pt depth 0pt\relax
            \vrule width 1pt height5.0pt depth 0pt\relax
            \vrule width 1pt height5.3pt depth 0pt\relax
            \vrule width 1pt height5.7pt depth 0pt\relax
            \vrule width 1pt height6.0pt depth 0pt\relax
            \vrule width 1pt height6.3pt depth 0pt\relax
            \vrule width 1pt height6.7pt depth 0pt\relax
            \vrule width 1pt height7.0pt depth 0pt\relax
            \vrule width 1pt height7.3pt depth 0pt\relax
            \vrule width 1pt height7.7pt depth 0pt\relax
            \vrule width 1pt height8.0pt depth 0pt\relax
            \vrule width 1pt height8.3pt depth 0pt\relax
            \vrule width 1pt height8.7pt depth 0pt\relax
            \vrule width 1pt height9pt depth 0pt\relax
      }


\def\widestnumber#1#2{}

\def\citewarning#1{\ifx\shlhetal\relax 
    \else
    \par{#1}\par
    \fi
}

\def\rm{\fam0 \tenrm}

\def\fakesubhead#1\endsubhead{\bigskip\noindent{\bf#1}\par}



%
%
%

%

\font\textrsfs=rsfs10
\font\scriptrsfs=rsfs7
\font\scriptscriptrsfs=rsfs5

\newfam\rsfsfam
\textfont\rsfsfam=\textrsfs
\scriptfont\rsfsfam=\scriptrsfs
\scriptscriptfont\rsfsfam=\scriptscriptrsfs

\edef\oldcatcodeofat{\the\catcode`\@}
\catcode`\@11

\def\Cal@@#1{\noaccents@ \fam \rsfsfam #1}

\catcode`\@\oldcatcodeofat


\expandafter\ifx \csname margininit\endcsname \relax\else\margininit\fi

\long\def\red#1\endred{}
\long\def\green#1\endgreen{}
\long\def\blue#1\endblue{}
\long\def\private#1\endprivate{}

\def\endred{ \unmatched endred! }
\def\endgreen{ \unmatched endgreen! }
\def\endblue{ \unmatched endblue! }
\def\endprivate{ \unmatched endprivate! }

\ifx\latexcolors\undefinedcs\def\latexcolors{}\fi

\def\emptycs{}
\def\evaluatelatexcolors{%
        \ifx\latexcolors\emptycs\else
        \expandafter\xxevaluate\latexcolors\xxfertig\evaluatelatexcolors\fi}
\def\xxevaluate#1,#2\xxfertig{\setupthiscolor{#1}%
        \def\latexcolors{#2}}


\font\smallfont=cmsl7
\def\rutgerscolor{\ifmmode\else\endgraf\fi\smallfont
\advance\leftskip0.5cm\relax}
\def\setupthiscolor#1{\edef\tmptmpcs{\noexpand\bgroup\noexpand\rutgerscolor
\noexpand\def\noexpand\currentcolor{#1}%
\noexpand}%
\expandafter\let\csname#1\endcsname\tmptmpcs
\def\tmptmpcs{\checkColorUnmatched{#1}\popthecolor}
\expandafter\let\csname end#1\endcsname\tmptmpcs}

\def\checkColorUnmatched#1{\def\expectcolor{#1}%
    \ifx\expectcolor\currentcolor   
    \else \edef\failhere{\noexpand\tryingToClose '\currentcolor' with end\expectcolor}\failhere\fi}

\def\currentcolor{???}

\def\popthecolor{\ifmmode\else\endgraf\fi\egroup}

\expandafter\def\csname#1\endcsname{}

\evaluatelatexcolors

 \let\outerhead\head
 \def\head{\innerhead}
 \let\innerhead\outerhead

 \let\outersubhead\subhead
 \def\subhead{\innersubhead}
 \let\innersubhead\outersubhead

 \let\outersubsubhead\subsubhead
 \def\subsubhead{\innersubsubhead}
 \let\innersubsubhead\outersubsubhead

 \let\outerproclaim\proclaim
 \def\proclaim{\innerproclaim}
 \let\innerproclaim\outerproclaim

 %
 %
 %
 %

\def\demo#1{\medskip\noindent{\it #1.\/}}
\def\enddemo{\smallskip}

\def\remark#1{\medskip\noindent{\it #1.\/}}
\def\endremark{\smallskip}

\pageheight{8.5truein}
\topmatter
\title{On $\lambda$ strongly homogeneity existence for cofinality logics} 
\endtitle
\author {Saharon Shelah \thanks {\null\newline I would like to thank 
Alice Leonhardt for the beautiful typing. Publication 750.
} \endthanks} \endauthor 

\affil{The Hebrew University of Jerusalem \\
Einstein Institute of Mathematics \\
Edmond J. Safra Campus, Givat Ram \\
Jerusalem 91904, Israel
 \medskip
 Department of Mathematics \\
 Hill Center-Busch Campus \\
  Rutgers, The State University of New Jersey \\
 110 Frelinghuysen Road \\
 Piscataway, NJ 08854-8019 USA} \endaffil
\medskip

\abstract  Let $C \underset\ne {}\to \subset \text{ Reg}$ be a
non-empty class (of regular cardinal).  Then the logic $\Bbb
L(Q^{\text{cf}}_C)$ has additional nice properties: it has
homogeneous model existence property.
\endabstract
\endtopmatter
\document

\newpage

\head {\S0 Introduction} \endhead  \resetall \sectno=0
 \spuriousreset
\bigskip

We deal with logics gotten by strengthening of first order logic by
generalized quantifiers, in particular, compact ones.
We continue \cite{Sh:199} (and \cite{Sh:43})

A natural quantifier is the cofinality quantifier, $Q^{\text{cf}}_{\le
\lambda}$ (or $Q^{\text{cf}}_C$), introduced in \cite{Sh:43} as the
first example of compact logic (stronger than first order logic, of
course). Recall that the ``uncountably many $x$'s" quantifier
$Q^{\text{card}}_{\ge \aleph_1}$, is $\aleph_0$-compact but not
compact.  But note that $\Bbb L(Q^{\text{cf}}_{\le \lambda})$ is a
very nice logic, e.g. with a nice axiomatization (in particular
finitely many schemes) like the one of $\Bbb L(Q^{\text{card}}_{\ge \aleph_0})$
of Keisler.  By \cite{Sh:199}, e.g. for $\lambda = 2^{\aleph_0}$,
its Beth closure is compact, giving the first compact logic with the
Beth property (i.e. implicit definition implies explicit definition).  

Earlier there were indications that having the Beth property is rare 
for such logic, see e.g. in Makowsky \cite{Mw85}.  A weaker version of
the Beth property is the weak Beth property dealing with implicit
definition which always works; H. Friedman claim that historically
this was the question.  Mekler-Shelah \cite{MkSh:166} prove that at
least consistently, $\Bbb L(Q^{\text{card}}_{\ge \aleph_1})$ satisfies
the weak Beth property.  V\"a\"an\"anen in the mid nineties motivated by the
result of Mekler-Shelah \cite{MkSh:166} asked whether we
can find a parallel proof for $\Bbb L(Q^{\text{cf}}_{\le \lambda})$.
\bigskip
\noindent
A natural property for logic ${\Cal L}$ is
\definition{\stag{0z.5} Definition}  A logic ${\Cal L}$ has the
homogeneous model existence property \ub{when} for every theory 
$T \subseteq {\Cal L}(\tau)$ (so has a model) has a 
strongly $({\Cal L},\aleph_0)$-homogeneous model $M$, i.e. $\tau_M =
\tau$ we have if $\bar a,\bar b \in {}^{\kappa >}M$ realizes 
the same ${\Cal L}(\tau)$-type
in $M$ then there is an automorphism of $M$ mapping $\bar a$ to $\bar b$.

It was introduced in 
\cite{Sh:199} as it helps to investigate the Beth property.

Adapting the proof of \cite{MkSh:166} was not enough.  Fine analysis
is needed.  So in \S1 we prove that $\Bbb L(Q^{\text{cf}}_C)$ has the
strongly $\aleph_0$-saturated model existence property.  
The situation concerning the weak Beth property is not clear.
\enddefinition
\bn
\margintag{0z.9}\ub{\stag{0z.9} Question}:  Does the logic $\Bbb L(Q^{\text{cf}(1)}_C)$
have the weak Beth property?

The first version of \S1 was done in 1996.
\bn
\margintag{0z.21}\ub{\stag{0z.21} Notation}:  1) $\tau$ denotes a vocabulary, ${\Cal L}$ a logic,
${\Cal L}(\tau)$ the language for the logic ${\Cal L}$ and the vocabulary
$\tau$.
\nl
2) Let $\Bbb L$ be first order logic, $\Bbb L(Q_*)$ be first order
   logic when we add the quantifier $Q_*$.
\nl
3) For a model $M$ and ultrafilter $D$ on a cardinal $\lambda$, let
   $M^\lambda/D$ be the ultrapower and $\bold j_{M,D} = \bold
   j^\lambda_{M,D}$ be the canonical embedding of $M$ into
   $M^\lambda/D$; of course, we can replace $\lambda$ by any set.
\nl
4) Let L.S.T. (theorem/argument) stand for L\"owenheim-Skolem-Tarski
   (on existence of elementary submodels).
\bn
Concerning \scite{0z.5}, more generally
\definition{\stag{0z.23} Definition}  1) $M$ is strongly 
$({\Cal L},\theta)$-saturated (in ${\Cal L} = \Bbb
   L$ we may write just $\theta$) \ub{when}
\mr
\item "{$(a)$}"  it is $\theta$-saturated (i.e. every set of ${\Cal
   L}(\tau_M)$-formulas with $< \theta$ parameters from $M$ and $<
   \theta$ free varialbes which is finitely satisfiable in $M$ and is
   realzied in $N$; if we omit it
\sn
\item "{$(b)$}"  if $\zeta < \theta$ and $\bar a,\bar b \in {}^\zeta
 M$ realizes the same ${\Cal L}(\tau_M)$-type in $N$, \ub{then} same
automorphisms of $M$ maps $\bar a$ to $\bar b$.
\ermn
2) $M$ is a strongly sequence $({\Cal L},\theta)$-homogeneous
   \ub{when} clause (b) above holds.
\nl
3) $M$ is sequence $(< \kappa)$-homogeneous when: if $\zeta <
   \kappa,\bar a \in {}^3 M,\bar b \in {}^3 M$ and tp$_\Delta(\bar
a,\emptyset,M) = \text{tp}_\Delta(\bar b,\emptyset,M)$ \ub{then} for
every $c \in M$ for some $d \in M$ we have tp$_\Delta(\bar b \char 94
\langle d \rangle,\emptyset,M) = \text{ tp}_\Delta(\bar a \char 94
\langle c \rangle,M)$.
\enddefinition
\bigskip

\definition{\stag{0z.26} Definition}  1) The logic $({\Cal L}$ has
``the $\kappa$-homogeneous existence property" \ub{when} every theory
$T \subseteq {\Cal L}(\tau_1)$ has a strongly $({\Cal
L},\kappa)$-homogeneous model.
\nl
2) Similarly ``the $\kappa$-saturated existence property.
\enddefinition
\bigskip

\newpage

\head {\S1 On strongly homogeneous models} \endhead  \resetall \sectno=1
 \spuriousreset
\bigskip

We prove that any theory in $\Bbb L(Q^{\text{cf}}_C)$ has strongly
$(\Bbb L(Q^{\text{cf}}_C),\theta)$-saturated models.
\bigskip

\definition{\stag{2b.1} Definition}   Let $\iota \in \{1,2\}$ and $C$
be a class of regular cardinals such that $C \ne \emptyset$, Reg.
\nl
1) The quantifier $Q^{\text{cf}(\iota)}_C$ is defined as follows:
\bn
\ub{syntactically}:  it bounds two variables, i.e. we can form
$(Q^{\text{cf}(\iota)}_C x,y)\varphi$, with its set of free variables
being defined as FVar$(\varphi) \backslash \{x,y\}$.
\bn
\ub{syntactically}:  $M \models (Q^{\text{cf}(\iota)}_C
x,y)\varphi(x,y,\bar a)$ \ub{iff} (a) + (b) holds where
\mr
\item "{$(a)$}"  relevancy demand: \ub{the case $\iota =1$}: the formula
$\varphi(-,-;\bar a)^M$ define so in $M$ a linear order with no last element
called $\le^\varphi_{M,\bar a}$ on the non-empty set
Dom$(\le^\varphi_{M,a}) = \{b \in M:M \models (\exists
y)(\varphi(b,y;\bar a)\}$
\sn
\ub{The case $\iota = 2$}: similarly but $\le^\varphi_{M,\bar a}$ is a
quasi linear order on its domain
\sn
\item "{$(b)$}"  the actual demand: $\le^\varphi_{M,\bar a}$ has
cofinality cf$(\le^\varphi_{M,\bar a})$, (necessarily an infinite
regular cardinal) which belongs to $C$.
\endroster
\enddefinition
\bn
\margintag{2b.1D}\ub{\stag{2b.1D} Convention}: 1) Writing 
$Q^{\text{cf}}_C$ we mean that this holds
for $Q^{\text{cf}(\iota)}_C$ for $\iota = 1$ and $\iota = 2$.
\nl
2) Let $\iota$-order mean order when $\iota = 1$ and quasi order when
   $\iota = 2$; but when we are using $Q^{\text{cf}}_C(\iota)$ then 
order means $\iota$-order.
\bigskip

\definition{\stag{2b.2} Definition}  1)  As $\{\psi \in 
{\Bbb L}(Q^{\text{cf}}_C):\psi$ has a model$\}$ does
not depend on $C$ (and is compact, see \cite{Sh:43}) we may use the
formal quantifier $Q_{\text{cf}}$, so the syntex is determined by not
the semantics, i.e. satisfaction relation $\models$.  We shall write 
$M \models_C \psi$ or $M \models_C T$ for the interpretation of 
$Q_{\text{cf}}$ as $Q^{\text{cf}}_C$, but can 
say ``$T \subseteq \Bbb L(Q_{\text{cf}})(\tau)$ has model/is
consistent".
\nl
2) If $C$ is clear from the context, \ub{then}
$Q^{\text{cf}}_\ell$ stands for $Q^{\text{cf}}_C$ if $\ell=1$ and
$Q^{\text{cf}}_{\text{Reg}\backslash C}$ if $\ell=0$.
\enddefinition
\bigskip

\demo{\stag{2b.4} Convention}  1) $T^*$ is a complete (consistent $\equiv$
has models) theory in ${\Bbb L}(Q_{\text{cf}})$ which is closed under
definitions i.e. every formula $\varphi = \varphi(\bar x)$ is equivalent to a
predicate $P_\varphi(\bar x)$.  
\nl
2) Let $T = T^* \cap$ first order logic, i.e. $T = T^* \cap 
\Bbb L(\tau_{T^*})$ is a complete first order theory.
\nl
3) $C \subseteq \text{ Reg}$, we let $C_1=C$ and 
$C_0 = \text{ Reg}\backslash C$, both non-empty.
\enddemo
\bigskip

\proclaim{\stag{2b.6} Theorem}  Assume $\chi = \text{\rm cf}(\chi),\mu
= \mu^{< \theta} \ge 2^{|T|} + \chi + \kappa,\theta \le \lambda,\theta
\le \text{ min}\{\chi,\kappa\},\chi \ne \kappa = \text{\rm cf}(\kappa)$ 
and
\nl
$\mu_\ell = \bigl\{\aligned \chi \quad
\ell=0 \\ \kappa \quad \ell=1 \endaligned$.

\ub{Then} there is a $\tau(T)$-model $M$ such that
\mr
\item "{$(a)$}"  $M \models T,\|M\| = \mu,M$ is $\theta$-saturated
\sn
\item "{$(b)$}"  if $\varphi(\bar z) = (Q^{\text{cf}}_\ell)
\psi(x,y;\bar z)$ then: $M \models P_{\varphi(\bar x)}(\bar a)$ iff
$\varphi(y,z;\bar a)$ define in $M$ a linear order with no last element
and cofinality $\lambda_\ell$
\sn
\item "{$(c)$}"   $M$ is strongly\footnote{as $T^*$ has elimination of
quantifiers, doing it for $\Bbb L(Q^{\text{cf}}_C)$ or for $\Bbb L$ is
 the same} $\theta$-saturated model of $T^*$.
\endroster
\endproclaim
\bigskip

\remark{Remark}  1) We can now change $\chi,\kappa,\mu$ and $\|M\|$ by
L.S.T.  Almost til the end instead $\mu \ge 2^{|T|} + \chi +
\kappa$ just $\mu \ge |T| + \chi + \kappa$ suffice.  The proof is
broken to a series of definitions and claims.  
The $\ge 2^{|T|}$ is necessary for $\aleph_0$-saturativity.
\nl
3) We can assume $\bold V$ satisfies GCH high enough and then use L.S.T.  So
$\mu^+ = 2^\mu$ below is not a real burden.
\endremark
\bigskip

\definition{\stag{2b.8} Definition}  0) Mod$_T$ is the class of models
of $T$.
\nl
1) 
\mr
\item "{$(a)$}"  $K = \{(M,N):M \prec N$ are from Mod$_T\}$
\sn
\item "{$(b)$}"  $K_\alpha = \{\bar M:\bar M = \langle M_i:i <
\alpha\rangle$ satisfies $M_i \in \text{ Mod}_T$ and
$i<j \Rightarrow M_i \prec M_j\}$ (so $K = K_2$)
\sn
\item "{$(c)$}"   $K^\alpha_\mu = \{\bar M \in K_\alpha:\|M_i\| 
\le \mu$ for $i < \alpha\}$, but then we (naturally) assume $\alpha < \mu^+$
\sn
\item "{$(d)$}"   let $\tau_\alpha = \tau_T \cup\{P_\beta:\beta <
\alpha\}$, each $P_\beta$ unary and $P_\alpha \notin \tau,\gamma <
\beta < \alpha \Rightarrow P_\gamma \ne P_\gamma\}$
\sn
\item "{$(e)$}"   for $\bar M \in K_\alpha$ let $\bold m(\bar M)$ be
the $\tau_\alpha$-model $M$ with universe $\cup\{M_\beta:\beta <
\alpha\},M \rest \tau_T = \cup\{M_\beta:\beta < \alpha\}$ and
$P^M_\beta = M_\beta$ and let $\bold m_0(\bar M)$ be the $\tau$-model
$\cup\{M_\beta:\beta < \alpha\}$ so $\bold m_0(\bar M) = \bold m(M)
\rest \tau$.
\ermn
2) Assume $(M^\ell,N^\ell) \in K$ for $\ell=1,2$ 
let $(M^1,N^1) \le (M^2,N^2)$ and
 $(M^1,N^1) \le_K (M^2,N^2)$ mean:
\mr
\item "{$(a)$}"  $M^1 \prec M^2$
\sn
\item "{$(b)$}"  $M^2 \cap N^1 = M^1$
\sn
\item "{$(c)$}"  $N^1 \prec N^2$
\sn
\item "{$(d)$}"  if $M^1 \models 
P_{(Q^{\text{cf}}_0,z)\varphi(y,z,\bar x)}[\bar a],c \in N^1,
c \in \text{ Dom}(\le^\varphi_{N_1,\bar a})$ and in
$N^1$ the element $c$ is $\le^\varphi_{M,\bar a}$-above all 
$d \in \text{ Dom}(\le^\varphi_{\bar M^1,\bar a})$, \ub{then} in
$N^2,c$ is $\le^\varphi_{N^2,\bar a}$-above all $d \in \text{ Dom}
(\le^\varphi_{M^2,\bar a})$.
\ermn
3) For $\bar M^1,\bar M^2 \in K_\alpha$ let $\bar M^1 \le \bar M^2$
   and $\bar M^1 \le_{K_\alpha} \bar M^2$ means $\gamma < \beta
\le \alpha \Rightarrow (M^1_\gamma,M^1_\beta) \le
   (M^2_\gamma,M^2_\beta)$.
\nl
4) For $\bar M \in K_\alpha,D$ an ultrafilter on $\lambda$ we define
$\bar N = \bar M^\lambda/D,\bold j_{M,D} = 
\bold j^\lambda_{\bar M,D}$ naturally: $N_\beta
 = M^\lambda_\beta/D$ for $\beta < \alpha$ and $\bold j_{\bar M,D} = 
\cup\{\bold j_{M_\beta,D}:\beta < \alpha\}$. 
\enddefinition
\bn
\margintag{2b.10}\ub{\stag{2b.10} Fact}:  0) For $\bar M^1,\bar M^2 \in K_\alpha$ we
have
\mr
\item "{$(a)$}"  $\bar M^1 \le_{K_\alpha} \bar M^2$ iff $\bold m(M^1)
\subseteq \bold m(\bar M^2)$
\sn
\item "{$(b)$}"  $(\bar M^\ell \rest P^{M_\beta}) \rest \tau_T =
M^\ell_\beta$.
\ermn
1) If $\bar M^1 \le_{K_\alpha} \bar M^2$ 
in $K_\alpha$ and $0 < \gamma < \beta \le \alpha$ \ub{then} 
$(\dbcu_{\varepsilon < \gamma} M^1_\varepsilon,\dbcu_{\varepsilon <
\beta} M^1_\varepsilon) \le (\dbcu_{\varepsilon < \gamma} M^2_\varepsilon,
\dbcu_{\varepsilon < \beta} M^2_\varepsilon)$ moreover
$\langle \dbcu_{i < 1+ \varepsilon} M^1_i:1+ \varepsilon 
\le \alpha\rangle \le_{K_{\alpha +1}} 
\langle \dbcu_{i < 1+ \varepsilon} M^2_i:1+\varepsilon \le \alpha\rangle$.
\nl
2) If $\langle \bar M^i:i < \delta\rangle$ is a
$\le_{K_\alpha}$-increasing sequence (of members of $K_\alpha$) and
we define $\bar M^\delta = \langle M^\delta_\varepsilon:\varepsilon <
\alpha\rangle$ by $M^\delta_\varepsilon = \cup\{M^i_\varepsilon:i <
\delta\}$ then $i < \delta \Rightarrow \bar M^i \le_{K_\alpha} \bar
M^\delta$ and the sequence 
$\langle \bar M^i:i \le \delta\rangle$ is continuous in $\delta$.
\nl
3) In part (2), if in addition $i < \delta \Rightarrow \bar M^i
   \le_{K_\alpha} \bar N$ so $\bar N \in K_\alpha$ \ub{then} $\bar
   M^\delta \le_{K_\alpha} \bar N$.
\nl
4) In part (2), if $\delta < \mu^+$ and $i < \delta \Rightarrow \bar
   M^i \in K^\alpha_\mu$ then $\bar M^\delta \in K^\alpha_\mu$.
\nl
5) If $\bar M \le_{K_\alpha} \bar N$ and $Y_\varepsilon \subseteq
N_\varepsilon$ for $\varepsilon < \alpha$ and
   $\Sigma\{\|M_\varepsilon\|,|Y_\varepsilon|:\varepsilon < \alpha\}
   \le \mu$ \ub{then} there is $\bar N' \in K^\alpha_\mu$ such that 
$\bar M \le_{K_\alpha} \bar N' \le \bar N$ 
 and $\varepsilon < \alpha \Rightarrow Y_\varepsilon \subseteq N'_\varepsilon$.
\nl
6) Assume $\bar M^i \in K^{\alpha(i)}_\mu$ for
$i < \delta_i < \mu^+,\langle \alpha(i):i < \delta\rangle$ is a
non-decreasing sequence of ordinals and $i < j < \delta \Rightarrow
\bar M^i \le_{K_{\alpha(i)}},\bar M^j \rest \alpha(i)$ and we define
$\alpha(\delta) = \cup\{\alpha(i):i < \delta\},\bar M^\delta = \langle
M^\delta_\beta:\beta < \alpha(\delta)\rangle$ where $M^\delta_\beta =
\cup\{M^i_\beta:\beta < \delta$ satisfies $\beta < \alpha(i)\}$ \ub{then}
$\bar M^\delta \in K^{\alpha(\delta)}_\mu$ and $i < \delta \Rightarrow
\bar M^i \le_{K_{\alpha(i)}} \bar M^\delta \rest \alpha(i)$.
\nl
7) If $\bar M^\ell \le_{K_\alpha} \bar N$ for $\ell =1,2$ and
$[a \in \bold m(\bar M^1) \Rightarrow a \in \bold m(\bar M^2)]$ then
in $\bar M^1 \le_{K_\alpha} \bar M^2$.
\bigskip

\demo{Proof}  Check.  \hfill$\square_{\scite{2b.10}}$
\enddemo
\bn
\margintag{2b.12}\ub{\stag{2b.12} Fact}:  1) If $(M_0,M_1) \in K_\mu$ and $(M_0,M'_1) \in K_\mu$
\ub{then} there are $M_2,f$ such that
\mr
\item "{$(a)$}"   $M'_1 \prec M_2 \in K_\mu$
\sn
\item "{$(b)$}"  $f \text{ is an elementary embedding of } 
M_1 \text{ into } M_2$
\sn
\item "{$(c)$}"  $f \rest M_0 = \text{ id}_{M_0}$
\sn
\item "{$(d)$}"  $(M_0,M'_1) \le_{K_2} (f(M_1),M_2)$.
\ermn
2) If $\bar M \in K_\alpha,\bar x = \langle x_\varepsilon:\varepsilon
< \zeta\rangle$ and $\Gamma$ is a set of first order formulas in
the variables $\bar x$ with parameters from the model $\bold m(\bar M)$ 
finitely satisfiable in $\bold m(M)$ such that $\varepsilon <
\zeta \Rightarrow \dsize \bigvee_{\beta < \alpha}
   P_\beta(x_\varepsilon) \in \Gamma$, \ub{then} there is $\bar N \in
   K_\alpha$ such that $\bar M \le_{K_\alpha} \bar N$ and $\Gamma$ is
   realized in $\bold m(\bar N)$.
\nl
3) If $\Gamma$ is a type over $\bold m_0(\bar M)$ of cardinality $<
   \text{ cf}(\alpha)$ \ub{then} it is included in some $\Gamma'$ as
   in part (2).
\nl
4) If $\bar M \in K^\alpha_\mu,D$ an ultrafilter on $\theta,M'_\beta =
M^\theta_\beta/D$, then $\bold j^\theta_{\bar M,D}(\bar M) 
\le_{K_\alpha} \langle M'_\beta:\beta < \alpha\rangle$, 
so for many $Y \in [\cup\{M'_\beta:\beta < \alpha\}]^\mu$
we have $\bold j^\theta_{\bar M,D}(\bar M) \le_{K_\alpha} 
\langle M'_\beta \rest Y:
\beta < \alpha\rangle \in K^\alpha_\mu$ \ub{but} see \scite{2b.12}.
\bigskip

\demo{Proof}  1) See \cite[\S4]{Sh:199}; just let $D$ be a 
regular ultrafilter on $\lambda \ge
\|M_1\|$, let $g$ an elementary embedding of $M_1$ into
$M^\lambda_0/D$ extending $\bold j = \bold j^\lambda_{M_0,D}$,
necessarily exists.

Lastly, let $M_2 \prec (M'_1)^\lambda/D$ include $\bold
j^\lambda_{M_1,D}(M'_1) \cup g(M_1)$ be of cardinality $\lambda$.
Identifying $M'_1$ with $\bold j^\lambda_{M'_1,D}(M'_1) < 
(M'_1)^\lambda/D$ we are done.
\nl
2) Similarly.
\nl
3) Trivial.
\nl
4) Should be clear.   \hfill$\square_{\scite{2b.12}}$ 
\enddemo
\bigskip

\definition{\stag{2b.14} Definition}  
$K^{\alpha,\text{ec}}_u$ is the class of $\bar M \in K^\alpha_\mu$
such that: if $\bar M \le_{K_\alpha} \bar N \in K^\alpha_\mu$, \ub{then}
$\bold m(\bar M) \le_{\Sigma_1} \bold m(\bar N)$, i.e.
\mr
\item "{$(*)$}"   if $a_1,\dotsc,a_n \in \bold m(\bar M),b_1,\dotsc,b_k \in
\bold m(\bar N),u \subseteq \alpha$ finite, $\varphi \in \Bbb L(\tau_T)$
is first order and $\bold m_0(\bar N) = 
\dbcu_{\beta < \alpha} N_\beta \models
\varphi(a_1,\dotsc,a_n,b_1,\dotsc,b_k)$ \ub{then} for some
$b'_1,\dotsc,b'_k \in \dbcu_{\beta < \alpha} M_\beta 
\text{ we have } \dbcu_\alpha M_\alpha \models
\varphi(a_1,\dotsc,a_n,b'_1,\dotsc,b'_k)$ and
$\dsize \bigwedge^k_{\ell=1} \,\, \dsize \bigwedge_{\beta \in  u} 
(b_\ell \in N_\beta \equiv b'_\ell \in M_\beta)$ 
(note that $a_m \in N_\beta \equiv a_m \in M_\beta$ by clause (b) of
Definition \scite{2b.8}(2).
\endroster
\enddefinition
\bigskip

\proclaim{\stag{2b.16} Claim}  1) $K^{\alpha,\text{ec}}_\mu$ is 
dense in $K^\alpha_\mu$ when $\mu \ge |\tau_T| + |\alpha|$. 
\nl
2) $K^{\alpha,\text{ec}}_\mu$ is closed under union of increasing
chains of length $< \mu^+$.
\endproclaim
\bigskip

\demo{Proof}  1) Given $\bar M_0 \in K^\alpha_\mu$ we try to choose $\bar
M_\varepsilon \in K^\alpha_\mu$ by induciton on $\varepsilon < mu^+$
such that $\langle \bar M_\zeta:\zeta \le \varepsilon\rangle$ is
$\le_{K_\alpha}$-increasing continuous and $\varepsilon = \zeta +1
\Rightarrow \bold m(\bar M_\zeta) \nleq_{\Sigma_1} \bold
m(M_\varepsilon)$.  For $\varepsilon = 0$ the sequence is given, for
$\varepsilon$ limit use \scite{2b.10}(2), for $\varepsilon = \zeta +1$
if we cannot choose then by \scite{2b.10}(5) we get $\bar M_\zeta \in
K^{\text{ec},\alpha}_M$ is as required.
\nl
2) Think on the definitions.  \hfill$\square_{\scite{2b.16}}$
\enddemo
\bigskip

\proclaim{\stag{2b.17} Claim}  1) If $\bar N \in K^{\text{ec}}_\alpha,Y
\subseteq \bold m_0(\bar N)$ and $\mu = |\tau_T| + |\alpha| + |Y|$
\ub{then} there is $\bar M \in K^{\text{ec},\alpha}_\mu$ such that
$\bar M \le_{K_\alpha} \bar N$.
\nl
2) If $\bar M^\ell \in K^\alpha_\mu$ and
$\bar M^0 \le_{K_\alpha} \bar M^\ell$ for $\ell=0,1,2$ and $\bar M^0
\in K^{\alpha,\text{ec}}_\mu$ \ub{then} we can find 
$(\bar N,f_1,f_2)$ such that:
\block
$\bar M^0 \le_{K_\alpha} \bar N\in K^\alpha_\mu$, moreover $\bar N \in
K^{\text{ec},\alpha}_\mu$ and $f_\ell$ is a
$\le_{K_\alpha}$-embedding of $\bar M^\ell$ into $\bar N$ over $\bar
M^0$.
\endblock
\endproclaim
\bigskip

\demo{Proof}  1) By the L.S.T. argument and part (2).
\nl
2) For $\ell=1,2$ let $\bar{\bold a}_\ell = \langle
a^\ell_\varepsilon:\varepsilon < \zeta_\ell\rangle$ list the elements
of $\bold m(\bar M^\ell)$ and let $\Gamma_\ell = \text{
tp}_{\text{qf}}(\bar{\bold a}_\ell,\emptyset,\bold m(\bar M^\ell) \equiv
\{\varphi(x^\ell_{\varepsilon_0},\dotsc,x^\ell_{\varepsilon_{n-1}},\bar
b):\varphi \in \Bbb L(\tau_T)$ is quantifier free, $\bar b \subseteq
\bold m(\bar M^0)$ and $\bold m(\bar M^\ell) \models
\varphi[a^\ell_{\varepsilon_0},\dotsc,a^\ell_{\varepsilon_{n-1}},
\bar b]\} \cup \{P_\beta(x^\ell_\varepsilon))^{\bold
t_\ell(\varepsilon,\beta)}:\beta < \alpha,\varepsilon < \zeta_\ell$
and $\bold t_\ell(\varepsilon,\beta)$ is the truth value of
$a^\ell_\varepsilon \in M^\ell_\beta\}$.

Now let $D$ be a regular ultrafilter on $\lambda = \|\bold m(\bar
M^1)\| + \|\bold m(\bar M^2)\|$ and use \scite{2b.12}(2),(3).  This is
fine to get $(f_1,f_2,\bar N)$ with $N \in K_\alpha$ and by
\scite{2b.10}(5) without loss of generality $N \in K^\alpha_\mu$ and by
\scite{2b.16}(1) \wilog \, $\bar N \in K^{\text{ec},\alpha}_\mu$.
\hfill$\square_{\scite{2b.17}}$ 
\enddemo
\bigskip

\proclaim{\stag{2b.18} Claim}  1) 
$(K^{\alpha,\text{ec}}_\mu,\le_{K^\alpha_\mu})$
has the JEP.
\nl
2) Suppose $\bar M^1,\bar M^2 \in K^\alpha_\mu,\beta \le \alpha,f$ 
is an elementary embedding of 
$\dbcu_{\gamma < \beta} M^1_\gamma$ into $\dbcu_{\gamma < \beta}
M^2_\gamma$ such that $\langle f(M_\gamma):\gamma < \beta\rangle  
\le_{K_\mu} \langle M^2_\gamma:\gamma < \beta\rangle$ (if   
$\beta=0$ then $f=\emptyset$ and there is no demand); equivalently $f$
is a $\le_{K_\beta}$-embedding of $\bold m(\bar M^1)$ into $\bold
m(\bar M^2 \rest B)$.

\ub{Then} we can find $M^3,f^+$ such that:
\mr
\item "{$(a)$}"  $\bar M^2 \le_{K_\mu} \bar M^3 \in K^\alpha_\mu$
\sn
\item "{$(b)$}"  $f \subseteq f^+$
\sn
\item "{$(c)$}"  $f^+ 
\text{ is an elementary embedding of } \dbcu_{\gamma < \alpha}
M^1_\gamma \text{ into } \dbcu_{\gamma < \alpha} M^3_\gamma$
\sn
\item "{$(d)$}"  $\langle f^+(M^1_\gamma)):\gamma < \alpha\rangle
\le_{K_\alpha} \langle M^3_\gamma:\gamma < \alpha\rangle$.
\endroster
\endproclaim
\bigskip

\demo{Proof}  1) A special case of part (2).
\nl
2) By induction on $\alpha$.
\mn
\ub{$\alpha =0$}:  nothing to do
\mn
\ub{$\beta = \alpha$}:  nothing to do
\mn
\ub{$\alpha =1$}:  so $\beta = 0$ which is trivial or $\beta =
\alpha$, a case done above.
\mn
\ub{$\alpha$ successor}:  by the induction hypothesis and transitive
nature of conclusion replacing $\bar M^2$ \wilog \, $\beta = \alpha
-1$, then use \scite{2b.12}.
\mn
\ub{$\alpha$ limit}:  By $\alpha - \beta$ successive uses of induction
hypothesis using \scite{2b.16}.  \hfill$\square_{\scite{2b.18}}$
\enddemo
\bigskip

\demo{\stag{2b.19} Conclusion}  $(K_\alpha,\le_{K_\alpha})$, or
formally ${\frak k} = (K_{\frak k},\le_{\frak k})$ defined 
by $K_{\frak k} := \{\bold m(\bar M):\bar M \in
K^{\text{ec}}_\alpha\},\bold m(\bar M^1) \le \bold m(\bar M^2)
\Leftrightarrow \bold m(M^1) \subseteq \bold m(\bar M^2)$, is an
   a.e.c. with amalgamation, and JEP, and LST$({\frak k}) \le |\tau_T|
   + |\alpha| + \aleph_0$.
\enddemo
\bigskip

\demo{Proof}  See \cite{Sh:88r} and history there.
\enddemo
\bn
\margintag{2b.20}\ub{\stag{2b.20} Fact}:  Assume $\lambda = \lambda^{< \lambda} >
|\tau_T| + \aleph_0 + |\alpha|$.  \ub{Then} there is $\bar M$ such
that
\mr
\item "{$(a)$}"  $\bar M \in K^{\text{ec}}_\alpha$ is universal for
 $(K^{\text{ec}}_\alpha,\le_{K_\alpha})$ in cardinality $\lambda$
\sn
\item "{$(b)$}"  $\bold m(\bar M)$ is model homogeneous for
 $(K^{\text{ec}}_\alpha,\le_{K_\alpha})$ of cardinality $\lambda$
\sn
\item "{$(c)$}"  $\bold m(\bar M)$ is sequence 
$(\Sigma_1,\lambda)$-homogeneous, see \scite{0z.23}(3).
\endroster
\bigskip

\demo{Proof}  Clause (a) + (b) are straight by \scite{2b.17} +
\scite{2b.18}(1), or use \scite{2b.18}(2), clause (c) follows: just
think.  \hfill$\square_{\scite{2b.20}}$.
\enddemo
\bn
\margintag{2b.22}\ub{\stag{2b.22} Fact}:  Assume $\bar M \in K^\alpha_\mu,\beta +1 <
\alpha,\ell \in \{0,1\}$ and $M_\beta \models 
P_{(Q^{\text{cf}}_\ell x,y)\varphi(x,y,\bar z)}[\bar a]$ 
\ub{then} there are $\bar N,c$ such that $\bar M \le \bar N \in
K^\alpha_\mu$ and:
\mr
\item "{$(*)_1$}"  if $\ell =1$ then $c \in 
\text{ Dom}(\le^\varphi_{N_\beta,\bar a})$ and $c$ is 
$\le^\varphi_{N_\gamma,\bar a}$-above 
$d \in \text{ Dom}(\le^\varphi_{M_\gamma,\bar a})$ for any  
$\gamma \in [\beta,\alpha)$
\sn
\item "{$(*)_2$}"  if $\ell =0$ then 
$c \in \text{ Dom}(\le^\varphi_{N_{\beta +1},\bar a})$ and 
is $\le^\varphi_{N_{\beta +1},\bar a}$-above any $d \in 
\text{ Dom}(\le^\varphi_{N_{\beta +1},\bar a})$.
\endroster
\bigskip

\demo{Proof}  First assume $\ell=1$, \wilog \, $\beta = 0$ as we 
can let $\bar N \rest \beta = M \rest \beta$.

By \scite{2b.18} \wilog \, $\alpha = 2$.  Now this is obvious by
\cite{Sh:43}, \cite{Sh:199}; in details by \cite{Sh:43} there is a
$\mu^+$-saturated model $M_*$ of $T$ such that $M_1 \prec M_*$ and
$M_* \models_{C_*} T^*$ whenever $\mu^{++} \in C_*$.  Let
$\{\varphi_i(x,y,\bar a^*_i):i < \mu\}$ list $\{\varphi(x,y,\bar
a'):\varphi \in \Bbb L(\tau_T),M_0 \models P_{\psi(Q_0^{c \ell}
x,y)\varphi(x,y;\bar z)}[\bar a']\}$, let $\langle
c_{i,\varepsilon},\varepsilon < \mu^+\rangle$ by
$\le^{\varphi_i,*}_{M_*,\bar a^*_i}$-increasing and cofinal.  For
$\varepsilon < \mu^+$ let $f_\varepsilon$ be an elementary embedding
of $M_1$ into $M_*$ over $M_0$ such that:
\mr
\item "{$(*)$}"  if $c \in \text{ Dom}(\le^{\varphi_i}_{M_*,\bar
a^*_i})$ is a $\le^{\varphi_i}_{M_*,\bar a^*_i})$-upper bound of
Dom$(\le^\varphi_i)_{M_0,\bar a_i})$ then $c_{i,\varepsilon}
\le^{\varphi_i}_{M_*,\bar a_i} c$.
\ermn
Let $c_* \in M_*$ be a $\le^\varphi_{M_*,\bar a}$-upper bound of
Dom$(\le^\varphi_{M_0,\bar a})$.  Choose $N_0 \prec M_*$ of cardinality
$\mu$ be such that $M_0 \cup\{c_*\} \subseteq N_0$ and choose
$\varepsilon < \mu^+$ large enough such that:
\mr
\item "{$(*)$}"  if $c < \mu$ and $d \in N_0$ is a $\le^{\varphi_i}_{M_*,\bar
a_i}$-upper bound of Dom$(\le^{\varphi_i}_{M_i,\bar a_i})$ then 
$\le^{\varphi_i}_{M_*,\bar a_i} c_{i,\varepsilon}$.
\ermn
Let $N_1 \prec M_*$ be of cardinality $\mu$ be such that $N_0 \cup
f_\varepsilon(M_1) \subseteq N_1$.  Renaming $f_\varepsilon$ is the
identity and $(N_0,N_1)$ is as required.

Second $\ell=0$ is even easier (again \wilog \, $\alpha =2$ and use
$N_0 = M_0,N_1$ realizes the relevant upper.  \hfill$\square_{\scite{2b.22}}$
\enddemo
\bigskip

\demo{\stag{2b.24} Conclusion}  In \scite{2b.20} the model 
$M^* = \bold m(\bar M^*) = \dbcu_{\beta < \alpha} M^*_\beta$ satisfies
\mr
\item "{$(a)$}"  if $M^* \models 
P_{(Q^{\text{cf}}_1 x,y)\varphi}[\bar a]$ then  the order
$\le^\varphi_{M^*,\bar a}$ has cofinality $\lambda$
\sn
\item "{$(b)$}"  if $\alpha$ is a limit ordinal and $M^* \models
P_{(Q^{\text{cf}}_0 x,y)\varphi}[\bar a]$ \ub{then} the linear order
$\le^\varphi_{M^*,\bar a}$ has cofinality cf$(\alpha)$
\sn
\item "{$(c)$}"  $M^*$ is cf$(\alpha)$-saturated
\sn
\item "{$(d)$}"  if $\mu^+ \in \bold C$ and cf$(\alpha) \in 
\text{ Reg}\backslash \bold C$ then $M^*$ is a model of $T^*$.
\endroster
\enddemo
\bigskip

\proclaim{\stag{2b.26} Claim}  Assume $\bar M \in
K^{\text{ec}}_\alpha$.
If $\zeta \le \mu$ and $\bar a,\bar b \in {}^\zeta(M^*_0)$ 
realize the same type (equivalently
q.f. type) in $M_0$ \ub{then} they realize the same $\Sigma_1$-type
in $\bold m(\bar M)$.
\endproclaim
\bigskip

\demo{Proof}  We choose $(N_\beta,f_\beta,g_\beta,h_\beta)$ by induction on
$\beta < \alpha$ such that:
\mr
\item "{$(a)$}"  $N_\beta$ is a model of $T$
\sn
\item "{$(b)$}"  $N_\beta$ is increasing continuous with $\beta$
\sn
\item "{$(c)$}"  $f_\beta,g_\beta$ are $\le_{K_{1+\beta}}$-embedding
of $\bar M \rest (1 + \beta)$ into $\langle N_\gamma:\gamma < 1 +
\beta\rangle \in K_{1 + \beta}$
\sn
\item "{$(d)$}"  $f_0(\bar a) = g_0(\bar b)$
\sn
\item "{$(e)$}"  if $\gamma < \beta$ then $f_\gamma \subseteq
f_\beta,g_\gamma \subseteq g_\beta$.
\ermn
For $\beta=0$ this speaks just on Mod$_T$.

For $\beta$ successor use \scite{2b.12}.

For $\beta$ limit as in the successor case, recalling we translated it
to the successor case (by \scite{2b.10}(1)).

Having carried the induction $f = \cup\{f_\beta:\beta < \alpha\}$ and
$g = \cup\{g_\beta:\beta < \alpha\}$ are $\le_{K_\alpha}$-embedding of
$\bar M$ into $\bar N = \langle N_\beta:\beta < \alpha\rangle$.  By
\scite{2b.16}(1) there is $\bar N' \in K^{\text{ec}}_\alpha$ which is
$\le_{K_\alpha}$-above $\bar N$.  Now as $\bar M \in
K^{\text{ec}}_\alpha$, the $\Sigma_1$-type of $\bar a$ in $\bold
m(\bar M)$ is equal to the $\Sigma_1$-type of $f(\bar a)$ in $\bold
m(\bar N')$, and the $\Sigma_1$-type of $\bar b$ in $\bold m(\bar
M)$ is equal to the $\Sigma_1$-type of $f(\bar a)$ in $\bold m(\bar
N')$.  But $f(\bar a) = f_0(\bar a) = g_0(\bar b) = g(\bar b)$, so we
have gotten the promised equality of $\Sigma_1$-types.
 \hfill$\square_{\scite{2b.26}}$
\enddemo
\bigskip

\demo{\stag{2b.27} Observation}  1) If $\bar M \in K^{\text{ec}}_\alpha$
and $\beta < \alpha$ \ub{then} $\bar M' := \bar M \rest [\beta,\alpha)
= \langle M_{\beta + \gamma}:\gamma < \alpha < \beta\rangle$ belongs
to $K^{\text{ec}}_{\alpha - \beta}$.
\nl
2) If $\bar M \in K_\alpha,\beta < \alpha$ and $\bar M \rest
   [\beta,\alpha) \le_{K_{\alpha,\beta}} \bar N'$ then for some $\bar
   N \in K_\alpha$ we have $\bar M \le_{K_\alpha} \bar N$ and $\bar N
   \rest [\beta,\alpha) = \bar N'$. 
\enddemo
\bigskip

\demo{Proof}  1) If not, then there is $\bar N' \in K_{\alpha - \beta}$
such that $\bar M' \le_{K_{\alpha-\beta}} \bar N'$ but $\bold m(\bar
M') \nleq_{\Sigma_1} \bold m(\bar N')$.  Define $\bar N = \langle
N_\gamma:\gamma < \alpha\rangle$ by: $N_\gamma$ is $M_\gamma$ if
$\gamma < \beta$ and is $N'_{\gamma-\beta}$ if $\gamma \in
[\beta,\alpha)$.  Easily $\bar M \le_{K_\alpha} \bar N \in K_\alpha$
but $\bold m(\bar M) \nleq_{\Sigma_1} \bold m(\bar N)$, contradiction
to the assumption $\bar M \in K^{\text{ec}}_\alpha$.
\nl
2) The proof is included in the proof of part (1).
\hfill$\square_{\scite{2b.26}}$ 
\enddemo
\bigskip

\proclaim{\stag{2b.28} Claim}  In \scite{2b.20} for each $\beta <
\alpha$ we have
\mr
\item "{$(a)$}"  $\langle M^*_{\beta + \gamma}:\gamma < \alpha -
\beta\rangle$ is homogeneous universal for $K^{\alpha-\beta}_\mu$
\sn
\item "{$(b)$}"  if $\alpha = \alpha - \beta$, i.e $\beta + \alpha =
\alpha$ \ub{then} there is an isomorphism
from $\bar M^*$ onto $\langle M^*_{\beta +\gamma}:\gamma < \alpha -
\beta\rangle$, in fact, we can determine $f(\bar a) = \bar b$ if $\bar
a \in {}^\zeta(M^*_0),\bar b \in {}^\zeta(M^*_\beta)$ and 
{\rm tp}$(\bar a,\emptyset,M^*_\beta) 
= \text{\rm tp}(\bar b,\emptyset,M^*_\beta)$.
\endroster
\endproclaim
\bigskip

\demo{Proof}  Chase arrows as usual recalling \scite{2b.27}.
\enddemo
\bn
\ub{Proof of Theorem \scite{2b.6}}:

Without loss of generality there is $\mu \ge \chi$ such that 
$2^\mu = \mu^+$ (why?  let $\mu > \chi$ be regular work
in $\bold V^{\text{Levy}(\mu^+,2^\mu)}$ and use
absoluteness argument, or choose set $A$ of ordinals such that 
${\Cal P}(\chi),T \in \bold L[A]$ and regular $\mu$ large enough such
that $\bold L[A] \models ``2^\mu = \mu^+"$, work in $\bold L[A]$ 
a little more) and for the desired conclusion (there is a model of 
cardinality $\chi$ such that ...) it makes no difference).
Let $\alpha = \kappa$ and let $\bar M^* \in
K^{\text{ec},\alpha}_\lambda$ be as in \scite{2b.10} for $\lambda :=
\mu^+$ and let $M_* = \cup\{M^*_\beta:\beta < \alpha\}$.  

Now
\mr
\item "{$(*)_1$}"  $M_*$ is a model of $T^*$ by the
$\{\mu^+\}$-interpretation.
\ermn
[Why?  By \scite{2b.24}.]
\mr
\item "{$(*)_2$}"  $M_*$ is $\kappa$-saturated.
\ermn
[Why?  Clearly $M^*_\beta$ is $\kappa$-saturated for each $\beta <
\kappa$.  As $\kappa$ is regular and $\langle M^*_\beta:\beta <
\kappa\rangle$ is increasing with union $M_*$, also $M_*$ is
$\kappa$-saturated.] 
\mr
\item "{$(*)_3$}"  $M_*$ is strongly $\aleph_0$-saturated and even
strongly $\kappa$-saturated, see Definition \scite{0z.23}(1).
\ermn
[Why?  Let $\zeta < \kappa$ and $\bar a,\bar b \in {}^\zeta(M_*)$ realize
the same q.f.-type (equivalent by first order type) in 
$M_*$.  As $\zeta < \kappa$
for some $\beta < \kappa$ we have $\bar a,\bar b \in
{}^\zeta(M_\beta)$.  Now by \scite{2b.28} we know that $\langle M^*_{\beta +
\gamma}:\gamma < \kappa\rangle \cong \langle M^*_\gamma:\gamma <
\kappa\rangle$, and by \scite{2b.26} the sequences $\bar a,\bar b$
realize the same $\Sigma_1$-type in $\bold m(\langle M^*_{\beta +
\gamma}:\gamma < \kappa\rangle)$ hence there is an automorphism $\pi$ of
it mapping $\bar a$ to $\bar b$.  So $\pi$ is an automorphism of $M_*$
mapping $\bar a$ to $\bar b$ as required.]

Lastly, we have to go back to models of cardinality $\mu = \mu^{<
\theta} \ge \lambda + \kappa + 2^{|T|}$, 
this is done by the L.S.T. argument recalling \scite{2b.24}.

More fully, let $\langle \bar M^\varepsilon:\varepsilon <
\lambda\rangle$ be $\le_{K^\alpha_\chi}$-increasing continuous
sequence with union $\bar M^*$.  For $\zeta < \theta$ and $\bar a,\bar
b \in {}^\zeta(M_*)$ let $f_{\bar a,\bar b}$ be an automorphism of
$M_*$ mapping $\bar a$ to $\bar b$.  Now the set of $\delta < \lambda$
satisfying $\circledast_\delta$ below is a club of $\lambda$ hence if
cf$(\delta) \ge \theta$ then $M = \cup\{M^\varepsilon_\beta:\beta <
\lambda\}$ is as required, where
\mr
\item "{$\circledast_\delta$}"  $(a) \quad$ if $\varepsilon <
\delta,\zeta < \theta$ and $\bar a,\bar b \in
{}^\zeta(\cup\{M^\zeta_\beta:\beta < \alpha\})$ realize the same
$\Sigma_1$-type in $\bar M^\zeta$ then $\cup\{M^\delta_\beta:\beta <
\alpha\}$ is closed under $f_{\bar a,\bar b}$ and under 
$f^{-1}_{\bar a,\bar b}$
\sn
\item "{${{}}$}"  $(b) \quad$ the witnesses for the cofinality work,
i.e.
{\roster
\itemitem{ $\bullet_1$ }  if $\beta < \alpha,\bar a \in {}^{\omega
>}(M^\delta_\beta),M^\delta_\beta \models P_{(Q^{\text{cf}}_0
 y,z)\varphi(y,z,\bar x)}[\bar a]$ then for some $\varepsilon <
 \delta$ we have $\bar a \subseteq M^\varepsilon_\beta$ and for every $\gamma
 \in (\beta,\alpha)$ there is $c = c_{\varphi,\bar a,\gamma} \in
 M^\varepsilon_{\gamma +1}$ which is a $\le^\varphi_{M^\varepsilon_{\gamma
 +1},\bar a}$-upper bound of
 Dom$(\le^\varphi_{M^\varepsilon_\gamma,\bar a})$, hence this holds
 for any $\varepsilon' \in [\varepsilon,\lambda)$
\sn
\itemitem{ $\bullet_2$ }  if $\beta < \alpha,\bar a \in 
{}^{\omega >}(M^\gamma_\beta)$ and $M^\delta_\beta \models 
P_{(Q^{\text{cf}}_1 y,z)\varphi(y,z,\bar x)}[\bar a]$ \ub{then} 
for arbitrarily large $\varepsilon <  \delta$ we have 
$\bar a \subseteq M^\varepsilon_\beta$ and there is
$c = c_{\varphi,\bar a} \in M^{\varepsilon +1}_\beta$ which is a
$(\le^\varphi_{M^{\varepsilon +1}_\gamma,\bar a})$-upper bound of
Dom$(\le^\varphi_{M^\varepsilon_\gamma})$ for every $\gamma \in
[\beta,\alpha)$.  \hfill$\square_{\scite{2b.6}}$ 
\endroster}
\endroster
\bigskip

\remark{Remark}  If you do not like the use of (set theoretic
absoluteness) you may do the following.  Use \scite{2b.35}, which is
legitimate as
\mr
\item "{$(a)$}"  the class $(K^{\text{ec}}_\alpha,\le_{K_\alpha})$ is
an a.e.c. with L.S.T. number $\le |T| + \aleph_0$ and amalgamation, so
\scite{2b.35}(1) apply
\sn
\item "{$(b)$}"  using $\Sigma_1$-types, it falls under \cite{Sh:3}
more exactly \cite{Sh:54}, so \scite{2b.35}(3) apply
\sn
\item "{$(c)$}"  we can define $K^{\text{ec}(\varepsilon)}_\alpha$ by
induction on $\varepsilon \le \omega$
\sn
\ub{$\varepsilon = 0$}:  $K_\alpha$
\sn
\ub{$\varepsilon = 1$}:  $K^{\text{ec}}_\alpha$
\sn
\ub{$\varepsilon = n+1$}:  $K^{\text{ec}(n+1)}_\alpha = \{\bar M \in
K^{\text{ec}(n)}_\alpha$: if $\bar M \subseteq N \in
K^{\text{ec}(n)}_\alpha$ then $\bold m(M) \le_{\Sigma_{n+1}} m(\bar N)\}$
\sn
\ub{$\varepsilon = \omega$}:  $K^{\text{ec}(\varepsilon)}_\alpha =
\cap\{K^{\text{ec}(n)}_\alpha:n < \omega\}$.
\ermn
On $K^{\text{ec}(\omega)}_\alpha$ apply \scite{2b.35}(2).
\endremark
\bigskip

\remark{\stag{2b.35} Remark}  1) Assume ${\frak k} = (K_{\frak
k},\le_{\frak k})$ is a a.e.c. with $\lambda >$ LST$({\frak k})$
and $\mu = \mu^{< \lambda}$.  For any $M \in K_\mu$ there is a
strongly model $\lambda$-homogeneous $N \in K_\mu$ which 
$\le_{\frak k}$-extend $M$, which means: if $M \in K_{\frak k}$ has
cardinality $< \lambda$ and $f_1,f_2$ are $\le_{\frak k}$-embedding of
$M$ into $N$ then for some automorphism $g$ of $N$ we have $f_2 = g
\circ f_1$.
\nl
2) Let $D$ be as in \cite{Sh:3} and $K_D =$ as below?  If $\lambda =
\lambda^{<\theta} \ge |D|$ and $M \in K_\lambda$ has cardinality
   $\lambda$ there is $N \in K_D$ of cardinality $\lambda$ which
   $\prec$-extend $M$ and is strongly $(D,\theta)$-homogenous, i.e.
\mr
 \item "{$(a)$}"  if $\zeta < \theta,\bar a,\bar b \in {}^\zeta N$
 realizes the same type \ub{then} some automorphism $f$ of $N$ maps
 $\bar a$ to $\bar b$
\sn 
\item "{$(b)$}"  $D = \{\text{\rm tp}(\bar a,\emptyset,N):\bar a
 \in {}^{\omega >} N\}$.
\ermn
3) Assume $\Delta \subseteq \Bbb L(\tau)$ not necessarily closed under
   negative $D$ is a set of $\Delta$-types, $K_D$ is class of
   $\tau$-model such that $\bar a \in {}^{\omega >} M \Rightarrow
   \text{\rm tp}(\bar a,\emptyset,M) \in D$ and $M \le_D N$ iff $M
   \subseteq N$ are from $K_D$ and $\bar a \in {}^{\omega >} M
   \Rightarrow \text{\rm tp}_\Delta(\bar a,\emptyset,M) = \text{\rm
   tp}_\Delta(\bar a,\emptyset,N)$.  Assume further $D$ is good,
   i.e. for every $M \in K_D$ and $\lambda$ there is a
   $(D,\lambda)$-sequence homogeneous model $N \in K_D$ which
   $\le_D$-extends $M$.  \ub{Then} for every $\lambda =  \lambda^{<
   \theta} > |T| + \aleph_0$ and 
$M \in K_D$ of cardinality $\lambda$ there is a strongly sequence
   $(\Delta,\lambda)$-homogeneous. 
\endremark
\bigskip

\demo{\stag{2b.33} Conclusion}  The logic $\Bbb L(Q^{\text{cf}}_C)$ has
   the homogeneous model existence property.
\enddemo
\bigskip

\demo{Proof}  Choose $\lambda \in C,\kappa \in \text{ Reg}\backslash
C$ and apply \scite{2b.6}.  \hfill$\square_{\scite{2b.33}}$
\enddemo

\nocite{ignore-this-bibtex-warning} 
\newpage
    
REFERENCES.  
\bibliographystyle{lit-plain}
\bibliography{lista,listb,listx,listf,liste}

\enddocument